\documentclass{amsart}
\usepackage{amsmath, amstext, amsxtra, amsthm, amscd, amsgen,amsbsy, amsopn, amsfonts, latexsym, amssymb, amscd}
\usepackage{marvosym}

\begin{document}



\def\wrt{with respect to}

\def\proofend{\hbox to 1em{\hss}\hfill $\blacksquare $\bigskip }

\newtheorem{theorem}{Theorem}[section]
\newtheorem{proposition}[theorem]{Proposition}
\newtheorem{lemma}[theorem]{Lemma}
\newtheorem{remark}[theorem]{Remark}
\newtheorem{remarks}[theorem]{Remarks}
\newtheorem{definition}[theorem]{Definition}
\newtheorem{corollary}[theorem]{Corollary}
\newtheorem{example}[theorem]{Example}
\newtheorem{assumption}[theorem]{Assumption}
\newtheorem{problem}[theorem]{Problem}
\newtheorem{question}[theorem]{Question}
\newtheorem{conjecture}[theorem]{Conjecture}
\newtheorem{rigiditytheorem}[theorem]{Rigidity Theorem}
\newtheorem{maintheorem}[theorem]{Main Theorem}
\newtheorem{mainlemma}[theorem]{Main Lemma}
\newtheorem{claim}[theorem]{Claim}

\newcommand{\sslash}{\mathbin{/\mkern-6mu/}}

\def\Z{{\mathbb Z}}
\def\R{{\mathbb R}}
\def\Q{{\mathbb Q}}
\def\C{{\mathbb C}}
\def\N{{\mathbb N}}
\def\H{{\mathbb H}}
\def\Zp #1{{\mathbb Z }/#1{\mathbb Z}}

\def\sec{{\rm sec}}
\def\diam{{\rm diam}}

\def\paperref#1#2#3#4#5#6{\text{#1:} #2, {\em #3} {\bf#4} (#5)#6}
\def\bookref#1#2#3#4#5#6{\text{#1:} {\em #2}, #3 #4 #5#6}
\def\preprintref#1#2#3#4{\text{#1:} #2 #3 (#4)}


\title[Nonnegative curvature and complex cohomology]{Nonnegative curvature, low cohomogeneity and complex cohomology}
\author{Anand Dessai}
\address{Department of Mathematics, Chemin du Mus\'ee 23, University of Fribourg, 1700 Fribourg (Switzerland) }
\email{anand.dessai@unifr.ch}
\urladdr{http://homeweb.unifr.ch/dessaia/pub/}
\date{\today}

\begin{abstract} We construct several infinite families of nonnegatively curved manifolds of low cohomogeneity and small dimension which can be distinguished by their cohomology rings. In particular, we exhibit an infinite family of eight-dimensional cohomogeneity one manifolds of nonnegative curvature with pairwise non-isomorphic complex cohomology rings.
\end{abstract}

\thanks{The work was partially supported by SNF Grant No. 200020-149761}
\maketitle

\noindent
\section{Introduction}\label{intro}

In this paper we give new information on the ``size" of the class of manifolds of nonnegative sectional curvature. Here the size will be measured in terms of the possible isomorphism types of cohomology rings. Our aim is to exhibit among these manifolds infinite families of small dimension and large symmetry which can be distinguished by their cohomology rings. In particular, we present in Theorem \ref{maintheorem} an infinite family of eight-dimensional cohomogeneity one manifolds of nonnegative curvature with pairwise non-isomorphic complex cohomology rings.
Throughout the paper we will restrict to closed simply connected manifolds. If not stated otherwise, curvature will refer to sectional curvature.

To begin with, let us briefly recall some existence and obstruction results for nonnegative curvature and a question of Grove which motivated our investigation.

Whereas only a few examples of manifolds with positive curvature are known many more nonnegatively curved examples have been constructed. This can be explained by the fact that certain constructions for nonnegative curvature do not hold, or are not known to hold, for positive curvature. In particular, the property of having nonnegative curvature is preserved under products and examples for nonnegatively curved manifolds are provided by all homogeneous spaces and biquotients, which are quotients of compact Lie groups. Moreover Grove and Ziller \cite{GZ00} have shown that among cohomogeneity one manifolds (i.e. manifolds with an action of
a Lie group with a codimension one orbit) there exist many examples which admit invariant metrics with nonnegative curvature. Despite the discrepancy between positively and nonnegatively curved examples, it is an open question whether there exist nonnegatively curved manifolds which do not admit a metric with positive curvature (recall that we restrict to simply connected manifolds). For a survey on constructions and examples we recommend \cite{Wi07, Zi07}.

A few obstructions to the existence of a nonnegatively curved metric are known. According to B\"ohm and Wilking \cite{BW07} any nonnegatively curved metric transforms under the Ricci flow to a metric of positive Ricci curvature, provided the fundamental group is finite. Hence, a nonnegatively curved manifold must satisfy the topological constrains imposed by positive Ricci and positive scalar curvature.

By Gromov's Betti number theorem \cite{Gr81} the sum of Betti numbers (\wrt \ any field of coefficients) of a Riemannian manifold is bounded from above by a constant depending only on the lower curvature bound, the upper diameter bound and the dimension. In particular, in any fixed dimension the sum of Betti numbers of nonnegatively curved Riemannian manifolds has a uniform upper bound. In other words the cohomology rings of such manifolds, viewed as graded vector spaces, belong to a finite number of isomorphism types and this number satisfies an upper bound which depends only on the dimension and is independent of the field of coefficients. The Betti number theorem gives a strong restriction on the class of manifolds of nonnegative curvature. A stronger restriction is implied by the so called Bott conjecture which states that any nonnegatively curved manifold is elliptic or at least rationally elliptic.

In \cite{Gr93} Grove asked whether in any fixed dimension the class of closed simply connected Riemannian manifolds satisfying uniform lower curvature and upper diameter bounds falls into only finitely many rational homotopy types. It follows from the Betti number theorem that this is the case in dimension $\leq 5$.

Grove's question has been answered into the negative first by Fang and Rong \cite{FR01} for lower negative curvature and upper diameter bounds and shortly after by Totaro \cite{To03} for nonnegatively curved manifolds. The examples of Fang and Rong are in any dimension $\geq 22$, satisfy uniform two-sided curvature bounds and can be distinguished already by their complex cohomology rings. Totaro's examples start in dimension $6$, which is the lowest possible dimension. His six-dimensional manifolds are nonnegatively curved biquotients with pairwise non-isomorphic rational cohomology rings (and, hence, are of different rational homotopy type). However, their real cohomology rings fall into only finitely many isomorphism types. Totaro also exhibits an infinite family in dimension $7$ with uniform two-sided curvature and upper diameter bounds and an infinite family in dimension $9$ with nonnegative curvature and uniform upper curvature and diameter bounds (see \cite{To03} for details). Again these manifolds can be distinguished by their rational cohomology rings, but their real cohomology rings fall into only finitely many isomorphism types.

In view of the examples above the size of the class of six-dimensional manifolds of nonnegative curvature is large with regard to their rational cohomology rings. The main purpose of this paper is to show that in slightly higher dimension this phenomenon already holds with regard to complex cohomology and under additional assumptions on the cohomogeneity. More precisely, we show

\begin{theorem}\label{maintheorem} There are infinitely many eight-dimensional simply connected Riemannian manifolds with nonnegative curvature, an isometric cohomogeneity one action and with pairwise non-isomorphic complex cohomology rings.
\end{theorem}

By taking products - for example with spheres - one gets the corresponding statement also in any dimension $\geq 10$.

We remark that the theorem above is sharp in several respects.
\begin{remarks}\label{firstintroremark}
\begin{enumerate}
\item From the classification of low dimensional simply connected homogeneous spaces (resp. cohomogeneity one manifolds) by Klaus \cite{Kl88} (resp. Hoelscher \cite{Ho10}) follows that the rational cohomology rings of simply connected homogeneous spaces (resp. cohomogeneity one manifolds) of dimension $\leq 8$ (resp. $\leq 7$) belong to only finitely many isomorphism types. Hence, the conclusion of the theorem fails in dimension $<8$ and fails for homogeneous spaces of dimension $\leq 8$.
\item The bound on the curvature in the theorem above cannot be changed from nonnegative to positive since Verdiani \cite{Ve04} has shown that an even-dimen\-sional manifold of positive curvature and isometric cohomogeneity one action is equivariantly diffeomorphic to a compact rank one symmetric space.
\end{enumerate}
\end{remarks}

The manifolds in Theorem \ref{maintheorem} are constructed as total spaces of $\C P^1$-bundles over the six-dimensional complex flag manifold. One can show that if one replaces in the construction the base space by any other homogeneous space of dimension $\leq 6$ then the real (and, hence, complex) cohomology rings of the total spaces fall into only finitely many isomorphism types.

It is not known (at least to the author) whether there exist infinite families of nonnegatively curved manifolds in dimension $6$ with pairwise non-isomorphic real or complex cohomology rings. In dimension $7$ the rational (resp. real) cohomology rings of simply connected rationally elliptic manifolds fall into infinitely (resp. only finitely) many isomorphism types (cf. \cite{He14I}). Hence, in view of the Bott conjecture one expects only finitely many real isomorphism types for seven-dimensional manifolds.

All manifolds in Theorem \ref{maintheorem} have second Betti number equal to $3$. It is not difficult to see that the construction does not lead to infinitely many real isomorphism types if the second Betti number of the total space is less than $3$ and the dimension is $\leq 8$. For smaller second Betti number we can show

\begin{theorem}\label{cohom two theorem}
In dimension $8$ (resp. $10$) there are infinitely many simply connected Riemannian manifolds with second Betti number equal to $2$, nonnegative sectional curvature, an isometric cohomogeneity two action and pairwise non-isomor\-phic rational (resp. complex) cohomology rings.
\end{theorem}

By taking products - for example with spheres - one obtains corresponding examples in higher dimensions. The manifolds in the theorem above are total spaces of $\C P^2$-bundles over complex projective spaces and examples of so-called generalized Bott manifolds, a special class of torus manifolds (a torus manifold is an orientable $2n$-dimensional manifold with an effective action by an $n$-dimensional torus with non-empty fixed point set). The isomorphism type of the integral cohomology ring of such manifolds has been studied extensively by Masuda and his coworkers in the context of cohomological rigidity problems (see for example \cite{CMS11}, \cite{CPS12}).

\begin{remark} The manifolds in Theorem \ref{cohom two theorem} can be described as quotients of a product of two spheres by free isometric torus-actions (see Proposition \ref{isoquotientprop}). This is no surprise since, according to recent work of Wiemeler (cf. \cite{Wi15}, Th. 1.2), any simply connected nonnegatively curved torus manifold is diffeomorphic to a quotient of a free linear torus action on a product of spheres.
\end{remark}

The manifolds in Theorem \ref{maintheorem} and Theorem \ref{cohom two theorem} all have positive Euler characteristic. By Cheeger's finiteness theorem \cite{Ch70} they do not admit metrics with uniform two-sided curvature and upper diameter bounds. We don't know whether there exist families of manifolds in these dimensions for which the conclusion in the theorems above still holds if one assumes in addition uniform upper curvature and diameter bounds.

In the theorems above the manifolds have small positive cohomogeneity. It would be interesting to determine the lowest possible dimension in which there are infinite homogeneous families (i.e. of cohomogeneity zero) with pairwise non-isomorphic cohomology rings (for coefficients $\Q$, $\R $ or $\C$). Recently, Herrmann \cite{He14II} has shown, among other things, that there are infinitely many simply connected $13$-dimensional homogeneous manifolds with pairwise non-isomorphic complex cohomology rings satisfying uniform upper curvature and diameter bounds.

\bigskip
The paper is structured as follows. In Section \ref{geometrysection} we describe an infinite family $\{M_{k,l}\}$  of eight-dimensional manifolds which is used in the proof of Theorem \ref{maintheorem}. The family consists of quotients of $SU(3)\times SU(2)$ by free isometric torus-actions. The geometrical and symmetry properties given in Theorem \ref{maintheorem} follow from this description. Section \ref{geometrysection} also contains a brief discussion of their symmetry rank. The manifolds $M_{k,l}$  can also be described as the total space of projective bundles associated to the sum of two complex line bundles over the complex flag manifold. In Sections \ref{topologysection} and \ref{section cohomologytheorem} we show that their complex cohomology rings represent infinitely many isomorphism types, thereby completing the proof of Theorem \ref{maintheorem}. In Section \ref{projective bundle section} we use $\C P^2$-bundles over complex projective spaces and some facts from number theory to prove Theorem \ref{cohom two theorem}.

\bigskip
\noindent
{\em Acknowledgements.} I would like to thank Wilderich Tuschmann for our numerous discussions on related questions in dimension less than $8$ which motivated this paper. Many thanks also to Martin Herrmann and Michael Wiemeler for interesting conversations about their recent work and for useful comments and to Fernando Galaz-Garcia for a helpful chat concerning manifolds with two-sided curvature bounds. I would also like to thank the referee for several suggestions which helped to improve the exposition of the paper.

\section{Geometric properties of the manifolds $M_{k,l}$}\label{geometrysection}

In this section we describe an infinite family $\{M_{k,l}\}$ of eight-dimensional manifolds used in the proof of Theorem \ref{maintheorem} and show the geometrical and symmetry properties stated. The manifolds $M_{k,l}$ are quotients of $SU(3)\times SU(2)$ by a free isometric action of a three-dimensional torus.

Let $T_{SU(3)}$ denote the standard maximal torus of $SU(3)$ given by unitary diagonal matrices of determinant one. We identify $T_{SU(3)}$ with the torus $T^2=S^1\times S^1$ via $T^2\to T_{SU(3)}$, $\mathrm{diag}(\lambda _1,\lambda_2)\mapsto \mathrm{diag}(\lambda _1,\lambda_2,\lambda_1^{-1}\cdot \lambda _2^{-1})$. Similarly, we identify the standard maximal torus $T_{SU(2)}$ of $SU(2)$ with $S^1$.

We equip $SU(3)$ and $SU(2)$ with bi-invariant Riemannian metrics. For $k,l\in \Z $, $(k,l)\neq (0,0)$, let $\rho_{k,l}$ be the homomorphism
$$\rho_{k,l}:T^2\to SU(2),\quad (\lambda_1,\lambda_2)\mapsto \mathrm{diag}(\lambda_1^{k}\cdot \lambda_2^{l}, \lambda_1^{-k}\cdot \lambda_2^{-l}).$$
Note that $\rho_{k,l}$ surjects onto $T_{SU(2)}\cong S^1$.

We next consider the action of the three-dimensional torus $T^3=T^2\times S^1\cong T_{SU(3)}\times T_{SU(2)}$ on $SU(3)\times SU(2)$ given by
\begin{small}$$T^3\times SU(3)\times SU(2)\to SU(3)\times SU(2),\; (t,s)(U_1,U_2):=(U_1\cdot t^{-1},\rho_{k,l}(t)\cdot U_2\cdot s^{-1}).$$\end{small}Note that $T^3$ acts freely and isometrically.

Let $M_{k,l}$ be the quotient manifold. From the construction we see that $M_{k,l}$ can be described as the total space of a bundle over the six-dimensional complex flag manifold $SU(3)/T_{SU(3)}$ with fiber $S^2$. The bundle is associated to the principal bundle $SU(3)\to SU(3)/T_{SU(3)}$ and the action of $T_{SU(3)}\cong T^2$ on $S^2=SU(2)/S^1$ induced by $\rho_{k,l}$.

We equip $M_{k,l}$ and $SU(3)/T_{SU(3)}$ with the submersion metrics. This gives the following sequence of Riemannian submersions
$$SU(3)\times SU(2)\overset{T^3}\to M_{k,l}\overset {S^2}\to SU(3)/T_{SU(3)}.$$

We will see below that the manifolds $M_{k,l}$ are of cohomogeneity one (i.e. are manifolds with an action of a Lie group $G$ with one dimensional orbit space). Let us recall that any
simply connected cohomogeneity one $G$-manifold admits a decomposition
$M=G\times _{K_-}D_{-}\, \cup \, G\times _{K_+}D_{+}$ as a union of two disk bundles, where $H\subset\{K_+, K_-\}\subset G$ are isotropy subgroups of $G$ and $D_{\pm}$ are disks with $\partial D_{\pm}=K_{\pm}/H$. Conversely, a group diagram $H\subset\{K_+, K_-\}\subset G$ where $K_{\pm}/H$ are spheres, defines a cohomogeneity one manifold (see for example \cite{GZ00} for details).
\bigskip

\begin{proposition}\label{geosummaryprop} $M_{k,l}$ is an eight-dimensional Riemannian manifold with nonnegative curvature and admits an isometric action by $SU(3)$ of cohomogeneity one.\end{proposition}

\noindent
{\bf Proof:} By the O'Neill formulas \cite{ON66} all spaces in the sequence of submersions above have nonnegative curvature. In addition the submersions are equivariant \wrt \ the isometric action of $SU(3)$ given by multiplication from the left.

Since $SU(3)$ acts transitively on $SU(3)/T_{SU(3)}$ and $\rho_{k,l}(T^2)$ acts on $S^2$ with one dimensional orbit space, we see that the action of $SU(3)$ on $M_{k,l}$ is of cohomogeneity one. More precisely, the isotropy groups are as follows.

Let $N:=\left (\begin{smallmatrix}1 & 0\\ 0 & 1\end{smallmatrix}\right )\cdot S^1$ and $S:=\left (\begin{smallmatrix}0 & 1\\ -1 & 0\end{smallmatrix}\right )\cdot S^1$ denote the fixed points of the action of $\rho_{k,l}(T^2)$ on $S^2$ and let $s_N$ and $s_S$ denote the corresponding sections in the bundle $M_{k,l}\to SU(3)/T_{SU(3)}$.

For a point which is in the image of the section $s_N$ or $s_S$ the isotropy is non-principal and conjugate to $T_{SU(3)}$ in $SU(3)$. Outside of the two sections the isotropy is principal and conjugate to $\rho_{k,l}^{-1}(\{\pm \mathrm{Id}\})\subset T^2\cong T_{SU(3)}$.\proofend

The manifolds $M_{k,l}$ can also be described as total spaces of projective bundles associated to a sum of two complex line bundles over the complex flag manifold $SU(3)/T_{SU(3)}$.

Let $S^1_i$ denote the $i$th factor in $T^2$, $i=1,2$, and let $\xi_i$ be the principal $S^1$-bundle over the complex flag manifold associated to the principal torus bundle $SU(3)\to SU(3)/T_{SU(3)}$ and the projection $T_{SU(3)}\cong T^2\to S^1_i$. In other words, $\xi _i$ is the principal $S^1$-bundle $SU(3)/S^1_j\to SU(3)/T_{SU(3)}$, $i\neq j$. Note that the principal torus bundle $SU(3)\to SU(3)/T_{SU(3)}$ is isomorphic to the sum of principal bundles $\xi _1\oplus \xi _2$.

Let ${\mathcal L}_i$ denote the complex line bundle associated to $\xi _i$. Consider the complex vector bundle $E:={\mathcal L}_1^{k}\otimes {\mathcal L}_2^{l}\oplus {\mathcal L}_1^{\hat k}\otimes {\mathcal L}_2^{\hat l}$ over $SU(3)/T_{SU(3)}$, $k,l,\hat k, \hat l\in \Z $.

By construction $E$ is isomorphic to the bundle $SU(3)\times _{\rho}\C ^2\to SU(3)/T_{SU(3)}$, where $T^2\cong T_{SU(3)}$ acts on $SU(3)$ by right multiplication and acts on $\C ^2$ via $\rho: T^2\to U(2)$, $\mathrm{diag}(\lambda _1,\lambda _2 )\mapsto \mathrm{diag}(\lambda_1^{k}\cdot \lambda_2^{l}, \lambda_1^{\hat k}\cdot \lambda_2^{\hat l})$. Passing to projective bundles we see that $P(E)$ is isomorphic to $SU(3)\times _{\rho}U(2)/T_{U(2)}\to SU(3)/T_{SU(3)}$.

Next suppose that $\hat k=- k$ and $\hat l=-l$. In this situation $\rho $ takes values in $SU(2)$ and is equal to $\rho_{k,l}$. The total space of $P(E)$ is isomorphic to $SU(3)\times _{\rho_{k,l}}SU(2)/S^1$ which is equal to $M_{k,l}$.
For latter reference we summarize the discussion in the following
\begin{lemma}\label{manifold is projective bundle lemma} Every $M_{k,l}$ can be described as the total space of a projective bundle associated to the sum of a complex line bundle over the complex flag manifold and its dual.\proofend
\end{lemma}

Using this description we will show in  the following two sections that the complex cohomology rings of the $M_{k,l}$ do not belong to finitely many isomorphism types.

We close this section with a brief discussion of the symmetry rank of the manifolds $M_{k,l}$. Recall that the symmetry rank of a Riemannian manifold is the rank of its isometry group \cite{GS94}.

Let us first note that the action of $T_{SU(3)}\times T_{SU(3)}$ on $SU(3)$ by left and right multiplication induces an ineffective isometric action on $M_{k,l}$ with one-dimensional kernel. This can be shown directly using the description of $M_{k,l}$ as total space of the fiber bundle $SU(3)\times _{\rho_{k,l}}SU(2)/S^1\to SU(3)/T_{SU(3)}$. Hence, the symmetry rank of $M_{k,l}$ is at least three. It follows from recent work of Wiemeler \cite{Wi15} that the symmetry rank cannot be larger for any Riemannian metric on $M_{k,l}$.

To explain this let us first note that $M_{k,l}$ is rationally elliptic since it is the quotient of $SU(3)\times SU(2)$ by a free torus action. Also, it follows from the description in Lemma \ref{manifold is projective bundle lemma} that the integral cohomology of $M_{k,l}$ vanishes in odd degrees. In particular, $M_{k,l}$ has positive Euler characteristic.

Suppose $M_{k,l}$ admits a smooth effective action by a four-dimensional torus. Since the Euler characteristic is non-zero the torus must act with fixed points. Thus, $M_{k,l}$ is a torus manifold which according to \cite{Wi15}, Th. 1.1, is homeomorphic to a quotient of a free linear torus action on a product of spheres. This implies that $SU(3)\times SU(2)$, the $2$-connected cover of $M_{k,l}$, is homeomorphic to a product of spheres contradicting the classical fact that $SU(3)$ is the total space of the {\em non-trivial} $S^3$-bundle over $S^5$.

\section{Cohomological properties of the manifolds $M_{k,l}$}\label{topologysection}
In this section we begin to investigate the cohomology of the manifolds $M_{k,l}$. Theorem \ref{cohomologytheorem} below rephrases the fact that their complex cohomology rings represent infinitely many isomorphism types. This gives the cohomological assertion of Theorem \ref{maintheorem}. This section contains some preliminary arguments. The proof of Theorem \ref{cohomologytheorem} will be completed in the following section.

Recall from Lemma \ref{manifold is projective bundle lemma} that the manifolds $M_{k,l}$ are total spaces of projective bundles associated to a sum of two complex line bundles over the complex flag manifold $SU(3)/T_{SU(3)}$. Their cohomology ring can be computed using the Leray-Hirsch theorem. We will review this in the general situation first and will specialize later to the projective bundles in question.

Let $\pi :E\to B$ be a complex vector bundle of rank $(r+1)$ over a manifold $B$ and let $P(E)$ be the projective bundle associated to $E$ (here and in the following we will allow ourselves to denote a bundle also by its total space). We also denote by $\pi$ the projection $P(E)\to B$.

Recall the following classical fact: If $L\to B$ is a complex line bundle then the projective bundles $P(E)$ and $P(E\otimes L)$ are canonically diffeomorphic. This follows directly using the description of vector bundles via cocycles, cf. for example \cite{GH94}, p. 515 (or by choosing a no-where vanishing, maybe non-continuous, section $\sigma :B\to L$ and by observing that the map $E\to E\otimes L$, $e\to e\otimes \sigma(\pi (e))$, defines a diffeomorphism $P(E)\to P(E\otimes L)$ which covers $\mathrm{id}_B$ and is independent of the choice of $\sigma$).

We denote by $y\in H^2(P(E);\Z)$ the negative of the first Chern class of the canonical line bundle over $P(E)$. By the Leray-Hirsch theorem $H^*(P(E);\Z )$ is a free $H^*(B;\Z )$-module (via $\pi^*$) with basis $(1,y,y^2,\ldots ,y^r)$. The cohomology ring $H^*(P(E);\Z )$ is isomorphic to (cf. for example \cite{GH94}, p. 606)
$$H^*(B;\Z )[y]/(y^{r+1}+c_1(E)\cdot y^{r}+c_2(E)\cdot y^{r-1}+\ldots +c_r(E)\cdot y+c_{r+1}(E)).$$

In the following we will assume that $E$ splits as a sum of a complex line bundle $L$ and a complex vector bundle of rank $r$. Then $\pi $ admits a section $s:B\to P(E)$ defined by mapping $b\in B$ to the fiber of $L\to B$ over $b$. Using the section we can split $H^*(P(E);\Z )$ as $\ker (s^*)\oplus \mathrm{im} (\pi ^*)\cong \ker (s^*)\oplus H^*(B;\Z )$.

As explained above $P(E)$ and $P(E\otimes L^{-1})$ are canonically diffeomorphic. In the cohomological computation for the projective bundles it will be convenient to replace $E$ by $E\otimes L^{-1}$. Doing so, we may assume that $E$ contains a trivial complex line bundle, denoted by $L_0$, as a summand. Note that in this situation $c_{r+1}(E)=0$ and $s^*:H^*(P(E);\Z )\to H^*(B;\Z )$ is induced by $y\mapsto 0$.

We now restrict to the situation where $E$ is the sum of the trivial line bundle $L_0$ and a line bundle $L_1$. Let $u:=c_1(E)=c_1(L_1)$ be the first Chern class and let $M_u$ be the total space of the associated $\C P^1$-bundle $\pi:P(E)\to B$. Note that the diffeomorphism type of $M_u$ is uniquely determined by the class $u\in H^2(B;\Z )$.

We will always identify $H^*(M_u;\Z)$ with $H^*(B;\Z)[y]/(y^2+u\cdot y)$ using the Leray-Hirsch theorem. More generally we will consider for any coefficient ring $R$ and any $u\in H^2(B;R)$ the graded ring $H^*_u:=H^*(B;R)[y]/(y^2+u\cdot y)$. Since we are interested in the isomorphism type of such rings let us record the following two elementary facts.

\begin{lemma}\label{unit lemma} Let $\lambda \in R^*$ be a unit. Then $H^*_u\cong H^*_{\lambda \cdot u}$.
\end{lemma}

\noindent
{\bf Proof:} Define $\Phi :H^*(B;R)[y]\to H^*(B;R)[y]$ by $y\mapsto \lambda ^{-1}\cdot y$ and as identity on $H^*(B;R)$. Then $\Phi(y^2+u\cdot y)=\lambda ^{-2}\cdot (y^2 +\lambda \cdot u\cdot y)$. Hence, $\Phi $ induces a well-defined isomorphism $H^*_u\to H^*_{\lambda \cdot u}$.\proofend

\noindent
We next note that a diffeomorphism $\phi:B\to B$ induces a bundle isomorphism $\phi ^*(E)\to E$ covering $\phi$, a diffeomorphism $M_{\phi^*(u)}\to M_u$ and an isomorphism $H^*_u\to H^*_{\phi^*(u)}$. Similarly one has

\begin{lemma}\label{automorphism lemma} Let $f$ be an automorphism of $H^*(B;R)$. Then $H^*_u\cong H^*_{f(u)}$.
\end{lemma}

\noindent
{\bf Proof:} Define $\Phi :H^*(B;R)[y]\to H^*(B;R)[y]$ by $y\mapsto y$ and as $f$ on $H^*(B;R)$. Then $\Phi(y^2+u\cdot y)=(y^2 +f(u)\cdot y)$. Hence, $\Phi $ induces a well-defined isomorphism $H^*_u\to H^*_{f(u)}$.\proofend

In the remaining part of this section we will assume that $B$ is the complex flag manifold $SU(3)/T_{SU(3)}$. The next lemma gives the connection to the manifolds $M_{k,l}$. Let ${\mathcal L}_1$ and ${\mathcal L}_2$ be the line bundles defined in the previous section and let $L_1:={\mathcal L}_1^{-2k}\otimes {\mathcal L}_2^{-2l}$.

\begin{lemma}\label{diffeo lemma} $M_{k,l}$ is diffeomorphic to $M_u$ for $u:=c_1(L_1)$.
\end{lemma}

\noindent
{\bf Proof:} Recall from the last section that $M_{k,l}\cong P({\mathcal L}_1^{k}\otimes {\mathcal L}_2^{l}\oplus {\mathcal L}_1^{-k}\otimes {\mathcal L}_2^{-l})$. Since the latter is diffeomorphic to the projective bundle associated to
$$L_0\oplus L_1\cong ({\mathcal L}_1^{k}\otimes {\mathcal L}_2^{l}\oplus {\mathcal L}_1^{-k}\otimes {\mathcal L}_2^{-l})\otimes ({\mathcal L}_1^{-k}\otimes {\mathcal L}_2^{-l})$$ it follows that $M_{k,l}$ and $P(L_0\oplus L_1)$ are diffeomorphic.\proofend

>From now on let $R=\C $. Thus, $H^*_u:=H^*(SU(3)/T_{SU(3)};\C)[y]/(y^2+u\cdot y)$ for $u\in H^2(SU(3)/T_{SU(3)};\C)$.

Let ${\mathcal P}:=P(H^2(SU(3)/T_{SU(3)};\C))$ be the space of complex lines in\linebreak $H^2(SU(3)/T_{SU(3)};\C)$. Note that ${\mathcal P}\cong \C P^1$ since $b_2(SU(3)/T_{SU(3)})=2$. By Lemma \ref{unit lemma} the isomorphism type of the ring $H^*_u$, $u\neq 0$, only depends on the line $\C \langle u\rangle \in {\mathcal P}$.

Two lines $\C \langle u\rangle ,\C \langle \widetilde u\rangle \in {\mathcal P}$ are called {\em equivalent} if $H^*_u\cong H^*_{\widetilde u}$. Thus, the isomorphism types of the rings $H^*_u$, $u\in H^2(SU(3)/T_{SU(3)};\C)$, $u\neq 0$, correspond to the equivalence classes in ${\mathcal P}$. If $u$ is an integral cohomology class we will call $\C \langle u\rangle$ an {\em integral} line and the equivalence class of $\C \langle u\rangle$ an {\em integral} equivalence class. We are now ready to state the main technical result of this paper.

\begin{theorem}\label{cohomologytheorem} Every integral equivalence class in ${\mathcal P}$ contains only finitely many integral lines.
\end{theorem}
\noindent
The proof will be given in the next section. Assuming this theorem we now prove Theorem \ref{maintheorem}.

\begin{theorem} There are infinitely many eight-dimensional simply connected Riemannian manifolds with nonnegative curvature, an isometric cohomogeneity one action and with pairwise non-isomorphic complex cohomology rings.

\end{theorem}

\bigskip
\noindent
{\bf Proof:} Consider the infinite family $\mathcal F$ of eight-manifolds $\{M_{k,l}\}$ where $k$ and $l$ are coprime positive integers.
By Proposition \ref{geosummaryprop} $M_{k,l}$ admits a Riemannian metric with nonnegative curvature and isometric action by $SU(3)$ of cohomogeneity one.

Recall from Lemma \ref{diffeo lemma} that $M_{k,l}$ is diffeomorphic to $M_u$ where $u:=c_1(L_1)$ and $L_1:={\mathcal L}_1^{-2k}\otimes {\mathcal L}_2^{-2l}$. Thus, $H^*(M_{k,l};\C)$ is isomorphic to $H^*_u$.

We note that $H^2(SU(3)/T_{SU(3)};\Z)$ is freely generated by $c_1({\mathcal L}_1)$ and $c_1({\mathcal L}_2)$.
Since $k$ and $l$ are coprime positive integers we see that different $u$ in this construction belong to different lines in $\mathcal P$.
According to Theorem \ref{cohomologytheorem} the equivalence class of $\C \langle u\rangle$ in ${\mathcal P}$ contains only finitely many integral lines. Hence, for fixed $u$ there are only finitely many manifolds in $\mathcal F$ with complex cohomology ring isomorphic to $H_u^*$.

Since $\mathcal F$ is infinite there exists an infinite subfamily with pairwise non-isomorphic complex cohomology rings.
\proofend

We would like to remark that the manifolds $M_{k,l}$ are closely related to the Aloff-Wallach spaces. The passage from one family to the other may be viewed as a sort of trade-off between good curvature/symmetry properties on the one side and richness of the cohomological type on the other. To explain this let us first recall that each $M_{k,l}$ is the total space of the $S^2$-bundle associated to a certain principal $S^1$-bundle over $SU(3)/T_{SU(3)}$ via the action of $S^1$ on $SU(2)/S^1\cong S^2$ induced by $\lambda \mapsto \mathrm{diag}(\lambda , \lambda ^{-1})$ (this follows from the description given in Section \ref{geometrysection}). The total spaces of the $S^1$-principal bundles all have isomorphic rational cohomology rings and admit, as shown by Aloff and Wallach \cite{AlWa75}, homogeneous metrics of positive curvature if $k\cdot l\cdot (k+l)\neq 0$. In contrast, the corresponding $M_{k,l}$ represent infinitely many non-isomorphic complex cohomology rings but have less good curvature/symmetry properties.

\bigskip

Before we begin with the proof of Theorem \ref{cohomologytheorem} we will first discuss some properties of the cohomology ring of $SU(3)/T_{SU(3)}$ and the action of the Weyl group.

We identify $SU(3)/T_{SU(3)}$ with $U(3)/T_{U(3)}$, where $T_{U(3)}$ denotes the standard maximal torus of $U(3)$ given by unitary diagonal matrices. Let us recall (cf. \cite{Bo53}) that $H^*(BT_{U(3)};\Z)\cong \Z[x_1,x_2,x_3]$, that the Weyl group $W$ of $U(3)$ acts on $H^*(BT_{U(3)};\Z)$ by permuting $x_1,x_2,x_3$ and that the integral cohomology of $U(3)/T_{U(3)}$ can be identified with the quotient of $H^*(BT_{U(3)};\Z)$ by the ideal generated by the Weyl-invariants of positive degree, i.e.
$$H^*(SU(3)/T_{SU(3)};\Z )\cong H^*(U(3)/T_{U(3)};\Z )\cong \Z [x_1,x_2,x_3]/(\sigma _1,\sigma _2,\sigma _3),$$ where $\sigma _i$ denotes the $i$th elementary symmetric function in $x_1,x_2,x_3$.

Hence, in terms of the basis $(x_1,x_2)$ of $H^2(SU(3)/T_{SU(3)};\Z )$ the integral cohomology ring of $SU(3)/T_{SU(3)}$ is isomorphic to
$$\Z [x_1,x_2]/(x_1^2+x_2^2+x_1\cdot x_2, x_1^2\cdot x_2+x_1\cdot x_2^2).$$ Note that $x_1^3, x_2^3$ and $x_1^2\cdot x_2^2$ belong to the ideal $(x_1^2+x_2^2+x_1\cdot x_2, x_1^2\cdot x_2+x_1\cdot x_2^2)$, and, hence, are zero in $H^*(SU(3)/T_{SU(3)};\Z )$.
In the following we will always identify  $H^*(SU(3)/T_{SU(3)};\C )$ with
$$\C [x_1,x_2]/(x_1^2+x_2^2+x_1\cdot x_2, x_1^2\cdot x_2+x_1\cdot x_2^2).$$

Any ring homomorphism $H^*(SU(3)/T_{SU(3)};\C )\to H^*(SU(3)/T_{SU(3)};\C )$ is determined by its restriction to $H^2(SU(3)/T_{SU(3)};\C )$ and we will use this linear map in the subsequent discussion. The latter will be described by a representing matrix $\left (\begin{smallmatrix} a_{11} & a_{12}\\a_{21} & a_{22}\end{smallmatrix}\right )$ for the basis $(x_1,x_2)$.

The action of the Weyl group $W$ on $H^*(BT_{U(3)};\Z)$ induces an action of $W$ on the cohomology rings $H^*(SU(3)/T_{SU(3)};\Z )$ and $H^*(SU(3)/T_{SU(3)};\C )$. Let us mention for completeness a more direct description of this action: The normalizer $N$ of $T_{SU(3)}$ in $SU(3)$ acts on $SU(3)/T_{SU(3)}$ via conjugation and this action induces the action of the Weyl group $W=N/T_{SU(3)}\cong S_3$ on the cohomology of $SU(3)/T_{SU(3)}$.

For latter reference we note that the action of the permutations
$$(1), (12), (13), (23), (123), (321)\in W$$ on $H^2(SU(3)/T_{SU(3)};\C )$ is represented by
$$\left (\begin{smallmatrix} 1& 0\\0 & 1\end{smallmatrix}\right ),\; \left (\begin{smallmatrix} 0& 1\\1 & 0\end{smallmatrix}\right ),\; \left (\begin{smallmatrix} -1& 0\\-1 & 1\end{smallmatrix}\right ),\; \left (\begin{smallmatrix} 1& -1\\0& -1\end{smallmatrix}\right ),\; \left (\begin{smallmatrix} 0 & -1\\1& -1\end{smallmatrix}\right ),\; \left (\begin{smallmatrix} -1 & 1\\-1 & 0\end{smallmatrix}\right ),$$
respectively.

The Weyl group acts by pre- and post-composition on the set of ring endomorphisms. The next two lemmas give representatives for the $W\times W$-orbits which will be important in the proof of Theorem \ref{cohomologytheorem}.

\begin{lemma}\label{basehomlemma1} Let $\Psi : H^*(SU(3)/T_{SU(3)};\C )\to H^*(SU(3)/T_{SU(3)};\C )$ be a ring homomorphism represented by $\left (\begin{smallmatrix} a_{11} & a_{12}\\a_{21} & a_{22}\end{smallmatrix}\right )$.

Then $a_{11}\cdot a_{21}\cdot (a_{11}-a_{21})=0$ and $a_{12}\cdot a_{22}\cdot (a_{12}-a_{22})=0$.
\end{lemma}

\noindent
{\bf Proof:} Consider the equations
$$\Psi (x_1^3)=(a_{11}\cdot x_1+a_{21}\cdot x_2)^3\text{ and }\Psi (x_2^3)=(a_{12}\cdot x_1+a_{22}\cdot x_2)^3.$$
Using $x_1^3=x_2^3=x_1^2\cdot x_2+x_1\cdot x_2^2=0$ it follows that $a_{11}\cdot a_{21}\cdot (a_{11}-a_{21})=0$ and $a_{12}\cdot a_{22}\cdot (a_{12}-a_{22})=0$.\proofend

\begin{lemma}\label{basehomlemma2} Let $\Psi : H^*(SU(3)/T_{SU(3)};\C )\to H^*(SU(3)/T_{SU(3)};\C )$ be a ring homomorphism. Then there exists $\omega _1,\omega _2\in W$ such that $\omega _1\circ \Psi \circ \omega _2$ is represented by either $\lambda \cdot \left (\begin{smallmatrix} 1& 0\\0 & 1\end{smallmatrix}\right )$, $\lambda \in \C ^*$, or by $\lambda \cdot \left (\begin{smallmatrix} 1& \varsigma \\ 0 & 0\end{smallmatrix}\right )$, $\lambda \in \C $, where $\varsigma$ satisfies $\varsigma^2+\varsigma +1=0$.
\end{lemma}

\noindent
{\bf Proof:}
Suppose $\Psi\neq 0$. It is easy to check that there exist $\tilde \omega _1,\tilde \omega _2\in W$ such that $\Psi_1:=\tilde \omega _1\circ \Psi \circ \tilde \omega _2$ is represented by a matrix $\left (\begin{smallmatrix} a_{11} & a_{12}\\a_{21} & a_{22}\end{smallmatrix}\right )$ with $a_{11}\neq 0$, $a_{12}\neq 0$ and $a_{21}\neq 0$.

Since $a_{11}\neq 0$ and $a_{21}\neq 0$ we get $a_{11}=a_{21}$ from the last lemma. Applying the lemma to $a_{12},a_{22}$ we see that $\Psi_1$ is either given by $\left (\begin{smallmatrix} a_{11} & a_{12}\\a_{11} & 0\end{smallmatrix}\right )$ or given by $\left (\begin{smallmatrix} a_{11} & a_{12}\\a_{11} & a_{12}\end{smallmatrix}\right )$.

In the first case $a_{12}=-a_{11}$ (apply the last lemma to $\left (\begin{smallmatrix} a_{11} & a_{12}\\a_{11} & 0\end{smallmatrix}\right )\cdot \left (\begin{smallmatrix} -1& 0\\-1 & 1\end{smallmatrix}\right )$). It follows that after composing $\Psi_1$ from the right with the element of $W$ corresponding to $\left (\begin{smallmatrix} 0 & -1\\1& -1\end{smallmatrix}\right )$ the homomorphism $\Psi _1$ transforms to a homomorphism represented by $\lambda \cdot \left (\begin{smallmatrix} 1& 0\\0 & 1\end{smallmatrix}\right )$, $\lambda \in \C ^*$.

In the second case $\Psi _1$ can be transformed (by composition of $\Psi_1$ from the left with the element corresponding to $\left (\begin{smallmatrix} 0 & -1\\1& -1\end{smallmatrix}\right )$) to a homomorphism $\Psi_2$ represented by $\left (\begin{smallmatrix} a& b\\0 & 0\end{smallmatrix}\right )$, for some $a,b\in \C ^*$. Using $\Psi_2(x_1)^2+\Psi_2(x_2)^2+\Psi_2(x_1)\cdot \Psi_2(x_2)=0$ one finds that $a^2+b^2+a\cdot b=0$.

Hence, up to the factor $a \in \C ^* $ the homomorphism $\Psi_2$ is represented by $\left (\begin{smallmatrix} 1& \varsigma \\ 0 & 0\end{smallmatrix}\right )$, where $\varsigma $ satisfies $\varsigma ^2+\varsigma +1=0$.\proofend

\section{Proof of Theorem \ref{cohomologytheorem}}\label{section cohomologytheorem}
We want to show that every integral equivalence class in ${\mathcal P}$ contains only finitely many integral lines. The idea of the proof is the following: Given non-zero integral classes $\widetilde u, u$ and an isomorphism $\Phi : H^*_{\widetilde u}\to H^*_u$ we will use the action of the Weyl group $W$ to change $\Phi$ into an isomorphism $\varphi : H^*_{\widetilde v}\to H^*_v$ which is in a suitable sense of standard form. Here $\widetilde v$ and $v$ are integral classes which are in the same $W$-orbit as $\widetilde u$ and $u$, respectively. We then show that $\C \langle v\rangle$ is determined by $\C \langle \widetilde v\rangle$ up to finite ambiguity. Since the Weyl group is finite we conclude that $\C \langle u\rangle$ is determined by $\C \langle \widetilde u\rangle$ up to finite ambiguity. Hence, the equivalence class of $\C \langle \widetilde u\rangle$ contains only finitely many integral lines.

Before we go into the proof let us recall from the last section: For any non-zero class $u\in H^2(SU(3)/T_{SU(3)};\C)$ the isomorphism type of $H^*_u$ only depends on the line $\C \langle u\rangle \in \mathcal P$ (see Lemma \ref{automorphism lemma}). The Weyl group $W$ acts on $H^2(SU(3)/T_{SU(3)};\C)$ and on $\mathcal P$. By Lemma \ref{diffeo lemma} an element $\omega \in W$ maps $H^*_u$ isomorphically to $H^*_{\omega(u)}$. In particular, two elements in $\mathcal P$ which belong to the same $W$-orbit are equivalent.

\bigskip
Let $\widetilde u, u \in H^2(SU(3)/T_{SU(3)};\Z)$ be non-zero classes such that $\C \langle \widetilde u\rangle$ and $\C \langle u\rangle$ are equivalent, i.e. $H^*_{\widetilde u}$ and $H^*_u$ are isomorphic. We will show that $\C \langle u\rangle$ is determined by $\C \langle \widetilde u\rangle$ up to finite ambiguity.

We denote by $\widetilde \pi ^*: H^*(SU(3)/T_{SU(3)};\C )\to H^*_{\widetilde u}$ the inclusion map and denote by $s^*: H^*_u\to H^*(SU(3)/T_{SU(3)};\C )$ the homomorphism induced by $y\mapsto 0$.

Let $\Phi :H^*_{\widetilde u}\to H^*_u$ be an isomorphism. Define
$$\Psi:=s^*\circ \Phi\circ \widetilde \pi^*:H^*(SU(3)/T_{SU(3)};\C )\to H^*(SU(3)/T_{SU(3)};\C ).$$ Since $\Phi$ is an isomorphism $\Psi$  does not vanish on $H^2(SU(3)/T_{SU(3)};\C )$. By Lemma \ref{basehomlemma2} there exist $\omega _1,\omega _2\in W$ and $\lambda \in \C ^*$ such that
$$\psi:=\omega _1\circ \Psi \circ \omega _2: H^*(SU(3)/T_{SU(3)};\C )\to H^*(SU(3)/T_{SU(3)};\C ) $$ is represented by either $\lambda \cdot \left (\begin{smallmatrix} 1& 0\\0 & 1\end{smallmatrix}\right )$ or by $\lambda \cdot \left (\begin{smallmatrix} 1& \varsigma \\ 0 & 0\end{smallmatrix}\right )$, where $\varsigma $ satisfies $\varsigma ^2+\varsigma +1=0$. After rescaling $\Phi $ (i.e. replace $\Phi(x)$, $x$ homogeneous, by $\lambda ^{-\deg(x )/2}\cdot \Phi(x)$) we can assume that $\lambda=1$.

Let $\widetilde v:=\omega _2^{-1}(\widetilde u)$, $v:=\omega _1(u)$ and let
$$\varphi:=\omega _1\circ \Phi \circ \omega _2:H^*_{\widetilde v}\to H^*_{\widetilde u}\to H^*_u\to H^*_v.$$
We note that $\psi$ is the homomorphism induced by $\varphi$. Note also that $\widetilde v$ and $v$ are integral cohomology classes.

The isomorphism $\varphi$ is determined by its restriction to $H^2_{\widetilde v}=\C \langle x_1,x_2,y\rangle $ which we represent by the matrix $A$ with respect to the basis $(x_1,x_2,y)$. From the discussion above $A$ takes the form $\left (\begin{smallmatrix} 1 & \varsigma & \alpha _1\\0 & 0 & \alpha_2 \\ b_1 & b_2 & \beta \end{smallmatrix}\right )$ or $\left (\begin{smallmatrix} 1 & 0 & \alpha _1\\0 & 1 & \alpha_2 \\ b_1 & b_2 & \beta \end{smallmatrix}\right )$. We will discuss the two cases separately.

\bigskip
\noindent
Let us first assume that $\varphi:H^*_{\widetilde v}\to H^*_v$ is represented by $A:=\left (\begin{smallmatrix} 1 & \varsigma & \alpha _1\\0 & 0 & \alpha_2 \\ b_1 & b_2 & \beta \end{smallmatrix}\right )$. Since $\varphi $ is an isomorphism $A$ is invertible. In particular, $(b_1,b_2)\neq (0,0)$.

\begin{lemma}\label{firstlemma} $\C \langle u\rangle$ is in the $W$-orbit of $\C\langle x_1\rangle$, $\C\langle x_2\rangle$, $\C\langle x_1-x_2\rangle$, $\C\langle x_1+2x_2\rangle$ or\linebreak $\C\langle 2x_1+x_2\rangle$.\end{lemma}

\noindent
{\bf Proof:} Suppose $\C \langle u\rangle$ is not in the $W$-orbit of $\C\langle x_1\rangle$, $\C\langle x_2\rangle$ and $\C\langle x_1+2x_2\rangle$. Note that  the same holds for $\C \langle v\rangle$ since $\C \langle v\rangle$ and $\C \langle u\rangle$ are in the same $W$-orbit. In particular, $v$ is a linear combination $\gamma _1\cdot x_1 + \gamma _2\cdot x_2$ with $\gamma _1$, $\gamma _2$ non-zero integers satisfying $2\gamma_1-\gamma_2\neq 0$.

Consider the relation
$$0=\varphi (x_1^2+x_2^2+x_1\cdot x_2)=(x_1+b_1\cdot y)^2 + (\varsigma \cdot x_1 +b_2\cdot y)^2 + (x_1+b_1\cdot y)\cdot (\varsigma \cdot x_1 +b_2\cdot y).$$ Using $y^2=-v\cdot y$, $\gamma _1\neq 0$ and $\varsigma ^2+\varsigma +1=0$ one finds
$$b_1^2+b_2^2 + b_1\cdot b_2=0,\quad b_2=-\frac {2+\varsigma}{1+2\varsigma}b_1 \quad \text{and}\quad b_i\neq 0.$$
Next consider the relation $0=\varphi (x_1)^3=(x_1+b_1\cdot y)^3$. Using $b_i\neq 0$, $\gamma_2\neq 0$, $2\gamma_1-\gamma_2\neq 0$ and $y^2=-v\cdot y$ one finds
$$b_1=\frac 3 {2\gamma _1-\gamma_2}\quad \text{and}\quad (\gamma_1+\gamma_2)\cdot (\gamma _1-2\gamma _2)=0.$$
Hence, $\C\langle v\rangle$ is equal to $\C\langle x_1-x_2\rangle$ or $\C\langle 2x_1+x_2\rangle$. This proves the lemma.\proofend

\bigskip
\noindent
Let us now assume that $\varphi:H^*_{\widetilde v}\to H^*_v$ is represented by $A:=\left (\begin{smallmatrix} 1 & 0 & \alpha _1\\0 & 1 & \alpha_2 \\ b_1 & b_2 & \beta \end{smallmatrix}\right )$.

\begin{lemma}\label{secondlemma} $\C \langle u\rangle$ is in the $W$-orbit of $\C \langle \widetilde u\rangle$, $\C\langle x_1\rangle$, $\C\langle x_1-x_2\rangle$, $\C\langle x_1+2x_2\rangle$, $\C\langle 2x_1+x_2\rangle$ or belongs to at most two other $W$-orbits.\end{lemma}

\noindent
{\bf Proof:} We first consider the relation $0=\varphi (y\cdot (y+\widetilde v))$. If we write $\varphi (y)=z+\beta \cdot y$, where $z:=\alpha _1\cdot x_1 + \alpha _2\cdot x_2$, and define $\gamma \in \C $ by $\varphi (\widetilde v)= \widetilde v + \gamma \cdot y$ then the relation is equivalent to
\begin{equation}\label{y^2 relation, part 1}z\cdot (z + \widetilde v)=0\end{equation} and
\begin{equation}\label{y^2 relation, part 2}\beta \cdot (\beta + \gamma)\cdot v=\beta \cdot \widetilde v + (2\beta +\gamma)\cdot z.\end{equation}
If $z=0$ we can conclude directly that $\beta \neq 0$ and $\C\langle v\rangle=\C\langle \widetilde v\rangle$.

If $z+\widetilde v=0$ we get $\beta \cdot (\beta + \gamma)\cdot v=-(\beta + \gamma )\cdot \widetilde v$. Note that $\beta + \gamma \neq 0$ since
$$0\neq \varphi(\widetilde v + y)=\widetilde v + \gamma \cdot y + z+\beta \cdot y=(\beta + \gamma )\cdot y.$$
Hence, $\C\langle v\rangle=\C\langle \widetilde v\rangle$.

In both cases we see that $\C \langle u\rangle$ is in the $W$-orbit of $\C \langle \widetilde u\rangle$.

\bigskip
Next consider the case that $z\neq 0$ and $z+\widetilde v \neq 0$. A computation (see Lemma \ref{zerodivisorlemma} below) shows that $z=\lambda _1\cdot x_{\pm}$ and $z+\widetilde v =\lambda _2\cdot x_{\mp}$, where $x_{\pm }:=x_1 +\frac 1 2 \cdot (1 \pm \sqrt {-3})\cdot x_2$ and $\lambda _i\in \C ^*$.

Using the relation
$$0=\varphi (x_1^2+x_2^2+x_1\cdot x_2)=(x_1+b_1\cdot y)^2 + (x_2 +b_2\cdot y)^2 + (x_1+b_1\cdot y)\cdot (x_2 +b_2\cdot y)$$
we find
\begin{equation}\label{bi equation} (b_1^2+b_2^2+b_1\cdot b_2)\cdot v=(2b_1 +b_2)\cdot x_1+(b_1+2b_2)\cdot x_2.\end{equation}

The remaining argument will be divided into two parts depending on whether $b_1^2+b_2^2+b_1\cdot b_2$ vanishes or not.

\bigskip
\noindent
{\em Claim:} If $b_1^2+b_2^2+b_1\cdot b_2=0$ then $\C\langle u\rangle$ belongs to at most two $W$-orbits.

\bigskip
{\em Proof:} Suppose $b_1^2+b_2^2+b_1\cdot b_2=0$. Then we have $b_1=b_2=0$. Hence, $\varphi (\widetilde v)=\widetilde v$ and $\gamma =0$. Also, $\beta \neq 0$, since $A$ is invertible.

Note that $\lambda_1$ and $\lambda _2$ are uniquely determined by $\widetilde v$ since $\widetilde v=\lambda _2\cdot x_{\mp}-\lambda _1\cdot x_{\pm}$ and $x_+,x_-$ form a basis of $H^2(SU(3)/T_{SU(3)};\C )$. Hence, $z$ is determined by $\widetilde v$ up to $\Z /2\Z$-ambiguity (more precisely, either $z=\lambda _1\cdot x_{+}$, where $\lambda_1$ is determined by $\widetilde v=\lambda _2\cdot x_{-}-\lambda _1\cdot x_{+}$ or $z=\lambda _1\cdot x_{-}$, where $\lambda_1$ is determined by $\widetilde v=\lambda _2\cdot x_{+}-\lambda _1\cdot x_{-}$).

Since $\gamma =0$ and $\beta\neq 0$ equation (\ref{y^2 relation, part 2}) gives $\beta\cdot v=\widetilde v + 2z$.
Since $z$ is determined by $\widetilde v$ up to $\Z /2\Z$-ambiguity we see that $\C\langle v\rangle$ is determined by $\C\langle \widetilde v\rangle$ up to $\Z /2\Z$-ambiguity. Hence, $\C\langle u\rangle$ belongs to at most two $W$-orbits in $\mathcal P$.\quad \checkmark

\bigskip
\noindent
{\em Claim:} If $b_1^2+b_2^2+b_1\cdot b_2\neq 0$ then $\C\langle u\rangle$ is in the $W$-orbit of  $\C\langle x_1\rangle$, $\C\langle x_1-x_2\rangle$, $\C\langle x_1+2x_2\rangle$ or $\C\langle 2x_1+x_2\rangle$.

\bigskip
{\em Proof:} Suppose $b_1^2+b_2^2+b_1\cdot b_2\neq 0$. By equation (\ref{bi equation}) we have
$$v=\frac{(2b_1 +b_2)\cdot x_1+(b_1+2b_2)\cdot x_2}{b_1^2+b_2^2+b_1\cdot b_2}.$$
Suppose $\C \langle u\rangle$ is not in the $W$-orbit of $\C\langle x_1\rangle$. Then the same holds for $\C \langle v\rangle$ since $\C \langle v\rangle$ and $\C \langle u\rangle$ are in the same $W$-orbit. In particular, $v$ is a linear combination $\gamma _1\cdot x_1 + \gamma _2\cdot x_2$ with $\gamma _1$, $\gamma _2\in \Z $ and $\gamma _2\neq 0$.

The relation $0=\varphi(x_1^3)$ is equivalent to
$$b_1\cdot (3x_1^2-3x_1\cdot b_1^2\cdot x_1\cdot v +b_1^2\cdot v^2)=0.$$
If $b_1=0$ then $b_2\neq 0$ and $v=\frac 1 {b_2}\cdot (x_1 +2x_2)$. Hence, $\C\langle v\rangle=\C\langle x_1+2x_2\rangle$.

Next assume $b_1\neq 0$. Recall that $v=\gamma_1\cdot x_1 + \gamma _2\cdot x_2$ and $\gamma _2\neq 0$. Suppose $\C\langle v\rangle\neq \C\langle x_1+2x_2\rangle$, i.e. $2\gamma_1-\gamma_2$. Using the same reasoning as in the proof of Lemma \ref{firstlemma} we conclude that the relation $0=\varphi(x_1^3)$ gives
$(\gamma_1+\gamma_2)\cdot (\gamma_1-2\gamma_2)=0$.
Hence, $\C\langle v\rangle=\C\langle x_1-x_2\rangle$ or $\C\langle v\rangle=\C\langle 2x_1+x_2\rangle$.

Thus, $\C\langle u\rangle$ is in the $W$-orbit of  $\C\langle x_1\rangle$, $\C\langle x_1-x_2\rangle$, $\C\langle x_1+2x_2\rangle$ or $\C\langle 2x_1+x_2\rangle$ as claimed.\quad \checkmark

\bigskip
\noindent
In view of the two claims above $\C\langle u\rangle$ is in the $W$-orbit of  $\C\langle x_1\rangle$, $\C\langle x_1-x_2\rangle$, $\C\langle x_1+2x_2\rangle$, $\C\langle 2x_1+x_2\rangle$ or belongs to at most two other $W$-orbits if $z\neq 0$ and $z+\widetilde v \neq 0$. This completes the proof of the lemma (modulo the proof of Lemma \ref{zerodivisorlemma}).\proofend

\bigskip
\noindent
In summary we have shown that if $H^*_{\widetilde u}$ and $H^*_u$ are isomorphic and $\widetilde u$ is fixed then $\C \langle u\rangle$ belongs to a finite number of $W$-orbits. Since the Weyl group is finite we conclude that $\C \langle u\rangle$ is determined by $\C \langle \widetilde u\rangle$ up to finite ambiguity. Hence, the equivalence class of $\C \langle \widetilde u\rangle$ contains only finitely many integral lines.

In order to complete the proof of Theorem \ref{cohomologytheorem} we are left to show the following

\begin{lemma}\label{zerodivisorlemma} Let $z_1, z_2\in H^2(SU(3)/T_{SU(3)};\C )$ be non-zero. If $z_1\cdot z_2=0$, then $z_1=\lambda_1\cdot x_{\pm}$ and $z_2=\lambda _2\cdot x_{\mp}$, for some $\lambda _1,\lambda_2\in \C ^*$, where $x_{\pm }:=x_1 +\frac 1 2 \cdot (1 \pm \sqrt {-3})\cdot x_2$.
\end{lemma}

\noindent
{\bf Proof:} Let $z_i=:A_i\cdot x_1+B_i\cdot x_2$, $i=1,2$. Using $z_i\neq 0$ and $z_1\cdot z_2=0$ one finds $A_i,B_i\neq 0$ for $i=1,2$. Let $\tilde z_i:=z/A_i=:x_1+C_i\cdot x_2$. Then we have
$$\tilde z_1\cdot \tilde z_2=0\iff C_1\cdot C_2=C_1+C_2=1$$
$$\iff C_1=\frac 1 2 \cdot (1 \pm \sqrt {-3}),\quad  C_2=\frac 1 2 \cdot (1 \mp \sqrt {-3}).$$\proofend

\section{Projective bundles over projective space}\label{projective bundle section}
In this section we prove Theorem \ref{cohom two theorem}. The manifolds which we will use are projective bundles associated to a sum of complex line bundles over a complex projective space. We begin with a more general description of some of their geometric properties which might be of independent interest.
\begin{proposition}\label{isoquotientprop}
Let $E$ be a complex vector bundles over $\C P^m$ and let $M=P(E)$ be the total space of the associated projective bundle. Suppose $E$ splits as a sum of $r+1$ complex line bundles.

Then $M$ is given as a quotient of $S^{2r+1}\times S^{2m+1}$ by a free action of a two-dimensional torus $T^2$. Moreover $S^{2r+1}\times S^{2m+1}$ admits a metric of nonnegative curvature such that $T^2$ acts by isometries. The quotient $M$ equipped with the submersion metric has nonnegative curvature and carries an ineffective isometric action by $U(m+1)\times T^{r+1}$ of cohomogeneity $r$.
\end{proposition}

For the manifolds used in the proof of Theorem \ref{cohom two theorem} we will choose $r=2$. We remark that the description of $M$  as a quotient of $S^{2r+1}\times S^{2m+1}$ in the proposition above remains valid for any complex vector bundle $E$ over $\C P^m$ of rank $r+1$. As will be shown the splitting of $E$ as a sum of complex line bundles allows to exhibit an ineffective isometric action by $U(m+1)\times T^{r+1}$ on $M$ which is of cohomogeneity~$r$.

\bigskip
\noindent
{\bf Proof:} We consider the principal $T^{r+1}$-bundle $P\to \C P^m$ associated to the direct sum decomposition of $E$ into complex line bundles and identify $M$ with\linebreak $P\times _{T^{r+1}}U(r+1)/(U(r)\times U(1))$, where $T^{r+1}$ acts on $U(r+1)/(U(r)\times U(1))$ from the left via the inclusion $T^{r+1}\hookrightarrow U(r+1)$ of the standard maximal torus.

Note that the transitive $U(m+1)$-action on $\C P^m$ from the left lifts canonically to a left action on the Hopf line bundle over $\C P^m$ and its powers and, hence, to any principal $S^1$-bundle over $\C P^m$. Thus, the homogeneous $U(m+1)$-action on $\C P^m$ lifts to the principal $T^{r+1}$-bundle $P\to \C P^m$. The existence of such a lift can also be deduced from general lifting properties in principal torus bundles. However in our situation everything can be made completely explicit and geometric.

We note that the $U(m+1)$-action and the principal $T^{r+1}$-action combine to a homogeneous $U(m+1)\times T^{r+1}$-action on $P$.

Next consider the Hopf fibration $\pi:S^{2m+1}\to \C P^m$. We recall that $\pi $ is the quotient map \wrt \ the action of the center $S^1\subset U(m+1)$ and that $\pi $ is $U(m+1)$-equivariant.

Let $\widetilde P:=\pi ^*(P)$ be the total space of the pullback bundle. Then $\widetilde P\to S^{2m+1}$ is an $U(m+1)$-equivariant principal $T^{r+1}$-bundle. Moreover $\widetilde P$ is homogeneous \wrt \ the action of $U(m+1)\times T^{r+1}$. The bundle map $\widetilde P\to P$ is given by taking the quotient \wrt \ the action of the center of $U(m+1)$.

Note that $\widetilde P\to S^{2m+1}$ is trivial as a non-equivariant principal torus bundle, i.e. isomorphic to $S^{2m+1}\times T^{r+1}\to S^{2m+1}$, since $H^2(S^{2m+1};\Z ^{r+1})=0$.

Next consider the associated sphere bundle
$$S^{2r+1}\hookrightarrow \widetilde P\times _{T^{r+1}}U(r+1)/U(r)\to S^{2m+1},$$ where $T^{r+1}$ acts from the left on $U(r+1)/U(r)$ via the inclusion $T^{r+1}\hookrightarrow U(r+1)$ of the standard maximal torus. From the above we conclude that the total space $N:=\widetilde P\times _{T^{r+1}}U(r+1)/U(r)$ is non-equivariantly diffeomorphic to $S^{2m+1}\times S^{2r+1}$. By construction $N$ comes with a free $T^2$-action given by the action of the center of $U(m+1)\times U(r+1)$. The quotient is equal to $P\times _{T^{r+1}}U(r+1)/(U(r)\times U(1))$. Hence, $M$ is diffeomorphic to the quotient of $S^{2m+1}\times S^{2r+1}$ by a free $T^2$-action.

Finally we observe that $U(m+1)\times T^{r+1}$ still acts (ineffectively) on $M$ with cohomogeneity $r=\dim _\R \C P^r -r$.

Let us now come to the statement about the curvature. Recall that $\widetilde P$ is a homogeneous $U(m+1)\times T^{r+1}$-manifold. Hence, we can identify $\widetilde P$ equivariantly with a quotient of $U(m+1)\times T^{r+1}$ and can equip $\widetilde P$ with a homogeneous metric of nonnegative curvature (e.g. the metric induced from a bi-invariant metric for $U(m+1)\times T^{r+1}$).

Similarly we can choose a metric on $U(r+1)/U(r)$ with nonnegative curvature such that $T^{r+1}$ acts isometrically (e.g. take the round metric on $S^{2r+1}$\linebreak $\cong U(r+1)/U(r))$.

With this choices the quotients $N=\widetilde P\times _{T^{r+1}}U(r+1)/U(r)\cong S^{2m+1}\times S^{2r+1}$ and $M=P\times _{T^{r+1}}U(r+1)/(U(r)\times U(1))\cong N/T^2$  inherit a metric of nonnegative curvature by the formulas of O'Neill. Moreover the free $T^2$-action on $N$ and the cohomogeneity $r$ action on $M$ are by isometries.\proofend

The manifolds which we use in the proof of Theorem \ref{cohom two theorem} are $\C P^2$-bundles over $\C P^2$ resp. $\C P^3$ and are of the type considered in the previous proposition. Hence, these manifolds carry a metric of nonnegative curvature and an isometric action of cohomogeneity two. The cohomological statement given in Theorem \ref{cohom two theorem} follows from the next two propositions.

\begin{proposition}\label{dim 8 cohom two theorem} There exists an infinite family of complex vector bundles $E_k\to \C P^2$, where each $E_k$ is a sum of three complex line bundles, such that the eight-dimensional manifolds $M_k:=P(E_k)$ have pairwise non-isomorphic rational cohomology rings.
\end{proposition}

\noindent
{\bf Proof:} We will use the following classical facts from number theory: Any prime $p\equiv 1 \bmod 3$ is of the form $d^2-d\cdot e +e^2$ for some $d,e\in \Z$ (cf. \cite{HW08}, Th. 254 on page 287). We also note that by Dirichlet's theorem on arithmetic progressions (cf. \cite{Se73}, Chap. VI, Cor. on page 74) there are infinitely many prime numbers congruent to $1$ modulo 3.

We choose an infinite strictly increasing sequence $(p_k)_k$ among these primes. For each $k$ we fix integers $d_k,e_k$ satisfying $d_k^2-d_k\cdot e_k +e_k^2=p_k$.

Let $E_k$ be the sum of three complex line bundles $L_0$, $L_1$ and $L_2$ over $\C P^2$ where $L_0$ is the trivial line bundle and $L_1$ and $L_2$ have first Chern class equal to $d_k\cdot x$ and $e_k\cdot x$, respectively. Here $x$ denotes a fixed generator of $H^2(\C P^2;\Z )$. By the Leray-Hirsch theorem the integral cohomology of $M_k:=P(E_k)$ is given by
$H^*(M_k;\Z)\cong \Z[x,y]/(x^3,\prod_{i=0}^2(y+u_i))$, where $u_0:=0$, $u_1:=d_k\cdot x$ and $u_2:=e_k\cdot x$, i.e. $u_i=c_1(L_i)$.

We want to show that the $M_k$ have pairwise non-isomorphic rational cohomology rings. Let $M_{\widetilde k}$ be another manifold, $\widetilde k\neq k$, and let $d_{\widetilde k}$, $e_{\widetilde k}$, $\widetilde u_i$ denote the corresponding parameters.

Suppose $\Phi :H^*(M_{\widetilde k};\Q )\to H^*(M_{k};\Q )$ is a ring isomorphism. The restriction of $\Phi$ to $H^2(M_{\widetilde k};\Q )=\Q\langle x,y\rangle$ defines a ring isomomorphism $\hat \Phi : \Q[x,y]\to \Q[x,y]$ which maps the ideal $(x^3, \prod_{i=0}^2(y+\widetilde u_i))$ onto $(x^3, \prod_{i=0}^2(y+u_i))$ and induces $\Phi$. Note that the ideals are generated by homogeneous elements of (cohomological) degree $6$.

Let $a,b \in \Q$ be defined by $\Phi(x)=a \cdot x +b \cdot y$. Then $\hat \Phi (x^3)=(a \cdot x +b \cdot y)^3$ must be of the form $\lambda \cdot x^3+\mu \cdot \prod_{i=0}^2(y+u_i)$ for some rational numbers $\lambda, \mu$. This gives
$$a^3=\lambda,\quad b ^3=\mu,\quad 3a^2\cdot b =b ^3\cdot d_k\cdot e_k\quad \text{and}\quad 3a\cdot b ^2 =b ^3\cdot (d_k+ e_k).$$
If $b \neq 0$ then the last two equations imply $d_k=e_k=0$ which gives a contradiction since $d_k^2-d_k\cdot e_k +e_k^2=p_k$. Hence, $b =0$ and $\Phi(x)=a \cdot x$.

Let $\alpha,\beta \in \Q$ be defined by $\Phi(y)=\alpha \cdot x +\beta \cdot y$. Since $\Phi$ is an isomorphism and $\Phi(x)=a \cdot x$ we have $a, \beta\neq 0$. Let us write $\hat \Phi(\prod_{i=0}^2(y+\widetilde u_i))$ as a linear combination $\hat \lambda \cdot x^3+\hat \mu \cdot \prod_{i=0}^2(y+u_i)$ for some rational numbers $\hat \lambda,\hat \mu$. Then we obtain the following relations
$$\hat \lambda = \alpha ^3+a \cdot \alpha ^2\cdot (d_{\widetilde k} + e_{\widetilde k})+a^2\cdot \alpha \cdot (d_{\widetilde k}\cdot  e_{\widetilde k}) ,\quad \hat \mu=\beta^3,$$
\begin{equation}\label{xysquare} 3\alpha +a \cdot (d_{\widetilde k} + e_{\widetilde k})=\beta\cdot (d_k+ e_k)\end{equation}
and
\begin{equation}\label{xsquarey} 3\alpha^2+2a\cdot \alpha \cdot (d_{\widetilde k} + e_{\widetilde k})+a ^2\cdot (d_{\widetilde k}\cdot  e_{\widetilde k})=\beta^2\cdot (d_k\cdot e_k).\end{equation}
If we solve for $\alpha $ in equation (\ref{xysquare}) and insert the result into equation (\ref{xsquarey}) we obtain
$$a ^2\cdot (d_{\widetilde k}^2-d_{\widetilde k}\cdot e_{\widetilde k}+e_{\widetilde k}^2) =\beta ^2\cdot (d_k^2-d_k\cdot e_k +e_k^2).$$
Since $a , \beta \neq 0$ and $p_{\widetilde k}=d_{\widetilde k}^2-d_{\widetilde k}\cdot e_{\widetilde k}+ e_{\widetilde k}^2, p_k=d_k^2-d_k\cdot e_k +e_k^2$ are different primes we arrive at a contradiction. Hence, the rational cohomology rings of $M_{\widetilde k}$ and $M_k$ are not isomorphic.\proofend

We remark that the arguments in the proof can be used to show that the conclusion of the proposition above fails if one replaces rational coefficients by real coefficients.

Another eight-dimensional family with pairwise non-isomorphic rational cohomology rings can be obtained by crossing Totaro's six-dimensional manifolds \cite{To03} with $S^2$.
The six-dimensional manifolds are biquotients of the form $(S^3)^3\sslash (S^1)^3$ and come with a (visible) cohomogeneity three action. Crossing with $S^2$ one obtains nonnegatively curved eight-dimensional manifolds of cohomogeneity three. One can show that their real cohomology rings fall into only finitely many isomorphism types. It would be interesting to know whether these manifolds also admit a cohomogeneity two action.

\bigskip
\noindent
We now turn to the proof of the statement in Theorem \ref{cohom two theorem} concerning ten-dimension\-al manifolds. The manifolds which we use are total spaces of projective bundles associated to sums of three complex line bundles over $\C P^3$. Their cohomology can be identified with the quotient of a polynomial algebra in two generators by an ideal generated by two homogeneous elements of {\em different} cohomological degree. This feature will simplify greatly the algebraic considerations.

\begin{proposition}\label{dim 10 cohom two theorem} There exists an infinite family of complex vector bundles $E_k\to \C P^3$, where each $E_k$ is a sum of three complex line bundles, such that the ten-dimensional manifolds $M_k:=P(E_k)$ have pairwise non-isomorphic complex cohomology rings.
\end{proposition}

\noindent
{\bf Proof:} Let $E$ be the sum of three complex line bundles $L_1$, $L_2$ and $L_3$ over $\C P^3$. Let $u_i:=c_1(L_i)$, $i=1,2,3$. By the Leray-Hirsch theorem the integral cohomology of $M:=P(E)$ is given by
$H^*(M;\Z)\cong \Z[x,y]/(x^4,\prod_{i=1}^3(y+u_i))$. Here $x$ denotes a generator of $H^2(\C P^3;\Z )$.

We want to show that the manifolds constructed in this way contain an infinite sequence $(M_k)_k$ with pairwise non-isomorphic complex cohomology rings.

Let $\widetilde M=P(\widetilde E)$ be another manifold and let $\widetilde u_i$, $i=1,2,3$, denote the corresponding first Chern classes.

Suppose $\Phi :H^*(\widetilde M;\C )\to H^*(M;\C )$ is an isomorphism of rings. The restriction of $\Phi$ to $H^2(\widetilde M;\C )=\C\langle x,y\rangle$ defines a ring isomomorphism $\hat \Phi : \C[x,y]\to \C[x,y]$ which maps the ideal $(x^4, \prod_{i=1}^3(y+\widetilde u_i))$ to $(x^4, \prod_{i=1}^3(y+u_i))$ and induces $\Phi$. Note that the ideals are generated by homogeneous elements of (cohomological) degree $8$ and $6$.

Hence the element $\prod_{i=1}^3(y+\widetilde u_i)$ (which is the one of smaller degree) must be mapped under $\hat \Phi $ to $C\cdot \prod_{i=1}^3(y+u_i)$, where $C\in \C $ is a constant.

Since the restriction of $\hat \Phi $ to $\C\langle x,y\rangle$ is an isomorphism and $\C[x,y]$ has no zero-divisors the constant $C\neq 0$.

We next note that by Gauss' lemma $\C[x,y]$ is a unique factorization domain and that the elements $(y+\widetilde u_i)$ and $(y+u_i)$ are irreducible.

Hence, after a permutation of the $u_i$ we may assume that $\hat \Phi ( y+\widetilde u_i)=C_i\cdot (y+u_i)$ for some $C_i\in \C ^*$. Let $\widetilde l_i,l_i\in \Z $ be defined by $\widetilde u_i=:\widetilde l_i\cdot x $ and $u_i=:l_i\cdot x $. If we write $\hat \Phi (x)=:a \cdot x +b \cdot y$ and $\hat \Phi(y)=:\alpha \cdot x +\beta \cdot y$ for complex numbers $a,b,\alpha,\beta $ we obtain the equations $(\alpha +\widetilde l_i\cdot a)=(\beta+\widetilde l_i\cdot b)\cdot l_i$ for $i=1,2,3$.

We claim that $b =0$ if the $l_i$ are pairwise different.
To see this consider $\hat \Phi (x^4)=(a \cdot x +b \cdot y)^4$ which belongs to the ideal $(x^4, \prod_{i=1}^3(y+u_i))$ and, hence, is of the form $(a \cdot x +b \cdot y)^4=\widetilde C\cdot x^4+g(x,y)\cdot \prod_{i=1}^3(y+u_i)$, where $\widetilde C\in \C $ is a constant and $g(x,y)\in \C[x,y]$ is homogeneous of degree $2$. If we specialize to $x=1$ we see that $\vert a +b \cdot y\vert $ is equal to $\vert \widetilde C\vert ^{1/4}$ for $y=-l_i$, $i=1,2,3$. In particular, $t\mapsto a + b \cdot t $ intersects $\{z\in\C \, \mid \, \vert z\vert =\vert \widetilde C\vert ^{1/4} \}$ for $t=l_i$, $i=1,2,3$. Since the $l_i$ are pairwise different it follows that $b =0$ (a line cannot intersect a circle in three different points).

So suppose $b=0$. Then $a,\beta\neq 0$, $l_i=\frac 1 \beta  \cdot (\alpha+\widetilde l_i\cdot a)$ and the defining parameters for $\widetilde M$ and $M$ are coupled by $l_i-l_j=\frac a \beta \cdot (\widetilde l_i - \widetilde l_j)$ for all $i,j$. It follows that there are infinitely many manifolds with pairwise non-isomorphic complex cohomology rings. A specific family is given by the manifolds $M_k$ which correspond to the parameters $\{l_1,l_2,l_3\}=\{0,1,k\}$, $k\geq 2$.\proofend


\bigskip

\begin{thebibliography}{9}
\bibitem[AlWa75]{AlWa75} S. Aloff and N. Wallach, {\em An infinite family of distinct 7-manifolds admitting positively
curved Riemannian structures}, Bull. Am. Math. Soc. 81 (1975) 93--97
\bibitem[BW07]{BW07} C. {B\"ohm} and B. {Wilking}, {\em Nonnegatively curved manifolds with finite fundamental groups admit metrics with positive Ricci curvature}, Geom. Funct. Anal. 17 (2007) 665--681
\bibitem[Bo53]{Bo53} A. Borel, {\em Sur la cohomologie des espaces fibr\'es principaux et des espaces homog\`enes de groupes de Lie compacts}, Ann. Math. 57 (1953) 115--207
\bibitem[Ch70]{Ch70} J. Cheeger, {\em Finiteness theorems for Riemannian manifolds}, Am. J. Math. 92 (1970) 61--74
\bibitem[CMS11]{CMS11} S. {Choi}, M. {Masuda} and D. Y. {Suh}, {\em Rigidity problems in toric topology: a survey}, Proc. Steklov Inst. Math. 275 (2011) 177--190
\bibitem[CPS12]{CPS12} S. Choi S. Park, D.Y. Suh, {\em Topological classification of quasitoric manifolds with the second Betti number 2}, Pac. J. Math. 256 (2012) 19--49
\bibitem[FR01]{FR01} F. Fang and X. Rong, {\em Curvature, diameter, homotopy groups, and cohomology rings,} Duke Math. J. 107 (2001) 135--158
\bibitem[GH94]{GH94} P. {Griffiths} and J. {Harris}, {\em Principles of algebraic geometry. 2nd ed.}, Wiley (1994)
\bibitem[Gr81]{Gr81} M. Gromov, {\em   Curvature, diameter and Betti numbers}, Comment. Math. Helv. 56 (1981) 179--195
\bibitem[Gr93]{Gr93}K. Grove, {\em Critical point theory for distance functions}, in: Differential geometry: Riemannian geometry, Proc. Symp. Pure Math. 54 (1993) 357--385
\bibitem[GS94]{GS94}K. Grove, C. Searle, {\em Positively curved manifolds with maximal symmetry-rank}, J. Pure Appl. Algebra 91 (1994) 137--142.
\bibitem[GZ00]{GZ00} K. Grove, W. Ziller: Curvature and symmetry of {M}ilnor
spheres. Ann. of Math. 152, (2000) 331-367
\bibitem[HW08]{HW08}G.H. {Hardy} and E.M. {Wright}, {\em An introduction to the theory of numbers, 6th ed.}, Oxford University Press (2008)
\bibitem[He14I]{He14I} M. Herrmann, {\em Classification and Characterization of rationally elliptic manifolds in low dimensions}, preprint, arXiv:1409.8036 (2014)
\bibitem[He14II]{He14II} M. Herrmann, {\em Homogeneous spaces, curvature and cohomology}, preprint,  arXiv:1411.1960v1 (2014)
\bibitem[Ho10]{Ho10} C.A. Hoelscher, {\em On the homology of low-dimensional cohomogeneity one manifolds}, Transform. Groups 15 (2010) 115--133
\bibitem[Kl88]{Kl88} St. Klaus, {\em Einfach zusammenh\"angende kompakte homogene R\"aume bis zur Dimension Neun},
Diploma Thesis, Univ. Mainz (1988) 98 p.
\bibitem[ON66]{ON66} B. {O'Neill}, {\em The fundamental equations of a submersion}, Mich. Math. J. 13 (1966) 459--469
\bibitem[Se73]{Se73} J-P. Serre, {\em A course in arithmetic}, Springer GTM 7 (1973)
\bibitem[To03]{To03} B. Totaro, {\em  Curvature, diameter, and quotient manifolds,} Math. Res. Lett. 10 (2003) 191--203
\bibitem[Ve04]{Ve04} L. {Verdiani}, {\em Cohomogeneity one manifolds of even dimension with strictly positive sectional curvature}, J. Differ. Geom. 68 (2004) 31--72
\bibitem[Wi15]{Wi15} M. Wiemeler, {\em Torus manifolds and non-negative curvature}, J. Lond. Math. Soc. 91 (2015) 667--692
\bibitem[Wi07]{Wi07} B. Wilking, {\em Nonnegatively and positively curved
manifolds}, in: Metric and comparison geometry. Surveys in Differential Geometry 11, International Press (2007) 25--62
\bibitem[Zi07]{Zi07} W. Ziller, {\em Examples of Riemannian manifolds with non-negative sectional curvature}, in: Metric and comparison geometry. Surveys in Differential Geometry 11, International Press (2007) 63--102
\end{thebibliography}
\end{document}